\documentclass{amsart}

\usepackage{amssymb,latexsym}
\usepackage{amsmath}
\usepackage{amsfonts}
\usepackage{amsthm}
\usepackage{amssymb}

\theoremstyle{plain}
\newtheorem{theorem}{Theorem}

\theoremstyle{definition}

\begin{document}

\title[Four squares from three numbers]{Four squares from three numbers}

\author[A. Dujella]{Andrej Dujella}
\address[Andrej Dujella]{
Department of Mathematics\\
Faculty of Science\\
University of \mbox{Zagreb}\\
Bijeni{\v c}ka cesta 30, 10000 Zagreb, Croatia
}

\author[L. Szalay]{L\'aszl\'o Szalay}
\address[L\'aszl\'o Szalay]{
Institute of Mathematics \\ University of Sopron \\
Bajcsy-Zsilinszky utca 4, Sopron, 9400, Hungary
}
\email{duje@math.hr}





\begin{abstract}
We show that there are infinitely many triples of positive integers $a,b,c$ (greater than $1$)
such that $ab+1$, $ac+1$, $bc+1$ and $abc+1$ are all perfect squares.
 \end{abstract}

\maketitle

\section{Introduction}

A \emph{Diophantine $m$-tuple} is a set of $m$ distinct positive integers
with the property that the product of any two of its distinct elements plus $1$ is a perfect square.
The first example of a Diophantine quadruple was found by Fermat, and it was the set $\{1,3,8,120\}$.
In 1969, Baker and Davenport \cite{baker-davenport} proved that Fermat's set cannot be extended to a Diophantine quintuple.
There are infinitely many Diophantine quadruples. E.g., $ \{k,k+2,4k+4,16k^3+48k^2+44k+12\}$
is a Diophantine quadruple for $k\geq 1$. In 2004, Dujella \cite{dujella-crelle}
proved that there is no Diophantine sextuple and that there are only finitely many Diophantine quintuples.
Finally, in 2019, He, Togb\'e and Ziegler \cite{he-togbe-ziegler}
proved that there is no Diophantine quintuple.

There are many known variants and generalizations of the notion of Diophantine $m$-tuples.
For a survey of various generalizations and the corresponding references see Section 1.5 of the book
\cite{dujella-springer}.

Here we will consider a variant that was introduced in several internet forums\footnotemark, %
\footnotetext[1]{\noindent{\tt https://www.mathpages.com/home/kmath481.htm}, \\
{\tt https://benvitalenum3ers.wordpress.com/2015/01/07/abc-ab1ac1bc1abc1-are-all-squares/}}%
and appeared also in Section 14.5 of the book \cite{brown-book}, where it is attributed to John Gowland.
We will consider triples of positive integers $a,b,c$ with the property that
$ab+1$, $ac+1$, $bc+1$ and $abc+1$ are perfect squares. Thus, we are interested in
Diophantine triples $\{a,b,c\}$ satisfying the additional property that $abc+1$ is also a perfect square.
If we allow that $a=1$, then the problem degenerates from four conditions to only three conditions that
$b+1$, $c+1$ and $bc+1$ are perfect squares, or in other words that $\{1,b,c\}$ is a Diophantine triple.
It is easy to see that there are infinitely many such triples, e.g. we may take $b=k^2-1$, $c=(k+1)^2-1$
for any $k\geq 2$. Hence, we will require that $a,b,c$ are positive integers greater than $1$.

Several examples of such triples were given in Section 14.5 of \cite{brown-book}, e.g.,
$(5,7,24)$, $(8,45,91)$, $(8,105,171)$, $(3,133,176)$, $(11,105,184)$, $(20,84,186)$,
$(44,102,280)$, \linebreak $(40,119,297)$, $(24,301,495)$, $(24,477,715)$.
However, it remained an open question whether there exist infinitely many such triples.
The main result of this paper gives an affirmative answer to that question.

\begin{theorem} \label{tm1}
There are infinitely many triples of positive integers $a,b,c$ greater than $1$
such that $ab+1$, $ac+1$, $bc+1$ and $abc+1$ are all perfect squares.
\end{theorem}

\section{The construction of infinitely many triples}

We will search for the solutions within so-called {\it regular Diophantine triples},
i.e., triples $\{a,b,c\}$, such that $c=a+b+2r$, where $ab+1=r^2$. Then $ac+1=(a+r)^2$ and $bc+1=(b+r)^2$,
so $\{a,b,c\}$ is indeed a Diophantine triple.
According to \cite{dujella-ramanujan}, most of Diophantine triples are of this form.

By studying and extending the list of known solutions, we can see that many of them
have the property that $a$ is of the form $A^2+4$:
\begin{align*}
(8&=2^2+4, 45, 91), \\
(8&=2^2+4, 105, 171), \\
(20&=4^2+4, 84, 186), \\
(40&=6^2+4, 119, 297), \\
(40&=6^2+4, 2387, 3045), \\
(85&=9^2+4, 672, 1235), \\
(85&=9^2+4, 11859, 13952), \\
(533&=23^2+4, 33475, 42456), \\
(533&=23^2+4, 509736, 543235), \\
(1160&=34^2+4, 165627, 194509), \\
(1160&=34^2+4, 2449135, 2556897), \\
(7400&=86^2+4, 7102165, 7568067), \\
(7400&=86^2+4, 101263737, 103002439), \\
(16133&=127^2+4, 34117191, 35617120), \\
(16133&=127^2+4, 482768440, 488366151).
\end{align*}
Almost all of these examples follow the following pattern:
$a$ is of the form $a=A_n^2+4$, where $A_n$ is a (two-sided) binary
recursive sequence defined by
$$ A_0=1, \quad A_1=6, \quad A_{n+1}=4A_{n}-A_{n-1}. $$
For $n\geq 1$, the elements of the sequence $A_n$ are: $6, 23, 86, 321, \ldots$,
while for $n\leq -1$, the elements of the sequence $-A_{-n}$ are: $2, 9, 34, 127, 474, \ldots$.

Next, we study the values of $r$ (from $ab+1=r^2$) in observed examples.
For each $a$, we had two triples with given property.
We will give details for the second (with larger $b$) triples.
We notice that $r$'s have the form
$r=A_n^2 R_n + A_{n+1} -2$, where
$$ R_0=2, \quad R_1=8, \quad R_n=4R_{n-1}-R_{n-2}+1, $$
and again we may extend the recurrence to negative indices,
so for $n\leq -1$, the elements of the sequence $R_{-n}$ are: $1,3,12,46,\ldots$.
(In the smaller triples, we have $r=A_n^2 R_{n-1} - A_{n-1} -2$.)

To simplify manipulations with the above introduced recursive sequences,
will we express them in the terms of the sequence
$$ P_0=0, \quad P_1=1, \quad P_n=4P_{n-1}-P_{n-2}. $$
The sequence $(P_n)$ satisfies
 \begin{equation} \label{binet}
P_n = \frac{1}{2\sqrt{3}}\left((2+\sqrt{3})^{n} - (2-\sqrt{3})^{n}\right).
\end{equation}
Let us denote $x=P_{n+1}$, $y=P_n$.
Then we have $A_n=x+2y$ and $R_n=\frac{1}{2}(5x-3y-1)$.

From (\ref{binet}), it follows easily that
\begin{equation} \label{hom}
x^2-4xy+y^2=1.
\end{equation}
We will use (\ref{hom}) to make further expressions as homogeneous as possible
in order to simplify expressions and in particular to allow factorizations.
In that way, we obtain
\begin{align*}
a&=5x^2-12xy+8y^2, \\
r&=\frac{17}{2} x^3-\frac{33}{2} x^2y-\frac{5}{2} x^2+14xy^2+6xy-7y^3-4y^2, \\[4pt]
b&=\frac{31}{2} x^4-\frac{55}{2} x^3y+\frac{75}{2} x^2y^2-25 xy^3+8 y^4-\frac{17}{2} x^3+\frac{33}{2} x^2y-14xy^2+7y^3, \\[4pt]
c&=\frac{31}{2} x^4-\frac{55}{2} x^3y+\frac{75}{2}x^2y^2-25xy^3+8y^4+\frac{17}{2} x^3-\frac{33}{2} x^2y+14xy^2-7y^3.
\end{align*}
In order to prove Theorem \ref{tm1}, it remains to check that $abc+1$ is a perfect square.
First we get
\begin{align*}
abc=&\frac{1}{4}(3y+8x)(5x^2-12xy+8y^2)(2y^3-2xy^2-2x^2y+3x^3) \\
&\mbox{}\times (10y^4-22xy^3+50x^2y^2-39x^3y+28x^4),
\end{align*}
and then by writing $1=(x^2-4xy+y^2)^5$ in $abc+1$, we finally obtain
\begin{equation*}
abc+1=\frac{1}{4}(22y^5-24xy^4-8x^2y^3+84x^3y^2-119x^4y+58x^5)^2,
\end{equation*}
which shows that $abc+1$ is indeed a perfect square.
\qed

\bigskip

For example, by taking $n=4$, we have $x=209$, $y=56$, and we get
$a=1435208$, $r=2347998213$, $b=3841321681771$, $c=3846019113405$,
and $abc+1=4604722693427179^2$.


\bigskip

{\bf Acknowledgements.} A.D. was supported by the Croatian Science Foundation under the project no. IP-2022-10-5008 and
the QuantiXLie Center of Excellence, a project co-financed by the
Croatian Government and European Union through the European Regional Development Fund (Grant PK.1.1.02.0004).


\begin{thebibliography}{99}

\bibitem{baker-davenport}
A. Baker and H. Davenport,
{\it The equations $3x^2 - 2 = y^2$ and $8x^2 - 7 = z^2$}, Quart. J. Math. Oxford Ser. (2) {\bf 20} (1969), 129--137.

\bibitem{brown-book}
K. Brown, Numbers, Author's edition, 2023.

\bibitem{dujella-crelle}
A. Dujella,
{\it There are only finitely many Diophantine quintuples}, J. Reine Angew. Math. {\bf 566} (2004), 183--214.

\bibitem{dujella-ramanujan}
A. Dujella, {\it On the number of Diophantine $m$-tuples}, Ramanujan J. {\bf 15} (2008), 37--46.

\bibitem{dujella-springer}
A. Dujella,
Diophantine $m$-tuples and Elliptic Curves, Springer, Cham, 2024.

\bibitem{he-togbe-ziegler}
B. He, A. Togb\'e and V. Ziegler,
{\it There is no Diophantine quintuple}, Trans. Amer. Math. Soc. {\bf 371} (2019), 6665--6709.


\end{thebibliography}
\end{document}